# REACHING THE BEST POSSIBLE RATE OF CONVERGENCE TO EQUILIBRIUM FOR SOLUTIONS OF KAC'S EQUATION VIA CENTRAL LIMIT THEOREM

By Emanuele Dolera, Ester Gabetta[1] and Eugenio Regazzini[2]

*Università degli Studi di Pavia, Università degli Studi di Pavia and Università degli Studi di Pavia*

Let $f(\cdot, t)$ be the probability density function which represents the solution of Kac's equation at time $t$, with initial data $f_0$, and let $g_\sigma$ be the Gaussian density with zero mean and variance $\sigma^2$, $\sigma^2$ being the value of the second moment of $f_0$. This is the first study which proves that the total variation distance between $f(\cdot, t)$ and $g_\sigma$ goes to zero, as $t \to +\infty$, with an exponential rate equal to $-1/4$. In the present paper, this fact is proved on the sole assumption that $f_0$ has finite fourth moment and its Fourier transform $\varphi_0$ satisfies $|\varphi_0(\xi)| = o(|\xi|^{-p})$ as $|\xi| \to +\infty$, for some $p > 0$. These hypotheses are definitely weaker than those considered so far in the state-of-the-art literature, which in any case, obtains less precise rates.

**1. Introduction.** This paper deals with the speed of approach to equilibrium—with respect to the variational distance—of the solution of the so-called *Kac equation*. Kac (1956) introduced a model, also known as Kac's caricature of a Maxwellian gas, for the motion of a single molecule in a chaotic bath of like molecules *moving on the line*. At time $t = 0$, the velocities are thought of as approximately stochastically independent and identically distributed with common law $\mu_0$ and finite second moment $m_2$. Kac formulated suitable hypotheses to prove that independence is preserved at later times. The Kac model provides the pattern for the analysis of certain more physically realistic kinetic models, first of all the Boltzmann equation for Maxwellian molecules. Many of the essential features of these more realistic models are, in any case, preserved in the Kac simplified setting. Accordingly, the specific

Received October 2007; revised April 2008.
[1]Supported by MUR Grant 06-015821 and by CNR-IMATI.
[2]Supported by MUR Grant 06-134525. Also affiliated with CNR-IMATI, Milano, Italy.

probabilistic methods we will use for the Kac model can be applied to the study of the asymptotic behavior of solutions of equations of Maxwellian molecules with constant collision kernel supported by a compact subset of $\mathbb{R}$. Moreover, the very same methods can have wide applicability to the approach to equilibrium of solutions of certain inelastic variants of the above-mentioned models, introduced in connection with the study of the behavior of granular materials and the redistribution of wealth in simple market economies. For background information on the physical aspects of Boltzmann's equation [see, e.g., Cercignani (1975) and Truesdell and Muncaster (1980)]. For inelastic variants of classical kinetics models, cf. Villani (2006), and for economic applications, see Matthes and Toscani (2008).

According to the Kac model, if $\mu_0$ is absolutely continuous with density, $f_0$, then the velocity of each particle at any time $t$ has a probability density function $f(\cdot, t)$, which is a solution of the *Boltzmann problem*

$$\frac{\partial f}{\partial t}(v,t) = \frac{1}{2\pi} \int_0^{2\pi} \int_{\mathbb{R}} [f(v\cos\theta - w\sin\theta, t) \\ \times f(v\sin\theta + w\cos\theta, t) - f(v,t) \cdot f(w,t)] \, dw \, d\theta \tag{1}$$

with initial data $f(\cdot, 0^+) = f_0(\cdot)$. Given any initial probability density function, the Cauchy problem (1) admits a unique solution. See, for example, Morgenstern (1954) and McKean (1966).

The present paper aims at providing sharp bounds for the variational distance

$$\|f(\cdot, t) - g_\sigma\|_1 := \int_{\mathbb{R}} |f(v,t) - g_\sigma(v)| \, dv \qquad (t \geq 0),$$

where $g_\sigma$ denotes the Gaussian probability density function having zero mean and variance $\sigma^2$, that is, $g_\sigma(v) = (\sigma\sqrt{2\pi})^{-1} \exp\{-v^2/2\sigma^2\}$ for $v$ in $\mathbb{R}$, and

$$\sigma^2 := \int_{\mathbb{R}} v^2 \, d\mu_0(v). \tag{2}$$

Its structure is as follows: Section 2 presents the main result along with a discussion about the assumptions made for its derivation. Section 3 contains preliminary notions and results to be used to prove the theorem formulated in Section 2. The main steps of the proof are explained in Section 4, while more technical details are deferred to the Appendix.

**2. Presentation of the main result.** Our main result is embodied in the following theorem.

THEOREM 2.1. *Assume that the initial probability density function, $f_0$, of Kac's equation (1) has finite fourth moment. Moreover, suppose*

$$(3) \qquad \varphi_0(\xi) := \int_{\mathbb{R}} e^{i\xi x} f_0(x)\,dx = o(|\xi|^{-p}) \qquad (|\xi| \to +\infty)$$

*for some strictly positive $p$. Then there is a constant $C$, depending only on the behavior of $f_0$, for which*

$$(4) \qquad \|f(\cdot,t) - g_\sigma\|_1 \leq C e^{-(1/4)t} \qquad (t \geq 0),$$

*where $\sigma$ is given as in (2).*

From Section 4 below, it follows that $C$ can be taken equal to

$2 + 2(\overline{n} + 2^{\overline{n}}\overline{n}!)$

$\qquad + \dfrac{1}{\sqrt{2}}\{\text{sum of the coefficients of } e^{-(1/4)t} \text{ in the right-hand sides}$

$\qquad\qquad \text{of } (26)\text{–}(29),\,(31)\} + 2\{\text{right-hand side of inequality } (36)\},$

where $\overline{n} := 9\lceil 2/p \rceil$ and $\lceil s \rceil$ stands for the least integer not less than $s$. Thus, $C$ depends only on the fourth moment $m_4$ of $f_0$, on $\sigma$, $p$ and $L_p := \sup_{\xi \in \mathbb{R}}\{|\xi|^p \cdot |\varphi_0(\xi)|\}$.

2.1. *Comparison with existing literature.* For the sake of expository completeness, it is worth situating the above theorem in existing literature. First, it should be recalled that *weak convergence* of the probability $\mu(\cdot, t)$ associated to $f(\cdot, t)$ holds true, as $t$ goes to infinity, if and only if $\mu_0$ has finite second moment [see Gabetta and Regazzini (2008)]. According to the Maxwell pioneering work, the limiting distribution is Gaussian $(0, \sigma^2)$ with $\sigma^2$ as in (2). Furthermore, when $m_2 = +\infty$, it has been proved that $\mu(\cdot, t)$ converges vaguely to the zero measure on $\mathbb{R}$ [see Carlen, Gabetta and Regazzini (2008)].

As to the rate of convergence when $m_2$ is finite, McKean gave significant results early on bounds for rates of convergence to equilibrium of solutions of (1), both in terms of Kolmogorov distance and in terms of variational distance [see McKean (1966)]. For recent developments concerning various weak metrics, see Gabetta and Regazzini (2006b). In any case, a relevant feature of McKean's way of reasoning is that it was based on central limit problem methods. In this paper, we use these very same methods also unlike what has happened in the last decades when authors have prevalently adopted strategies of an analytical nature. See, for example, the recent survey in Villani (2008).

Apropos of the variational distance, McKean conjectured that $\|f(\cdot, t) - g_\sigma\|_1$ has an upper bound which goes to zero exponentially with a rate determined by the *least negative eigenvalue* (equal to $-1/4$) of the linearized

collision operator. Compared to these statements, Theorem 2.1 represents the first satisfactory support of the McKean assertion. Indeed, the problem of validating this conjecture has been seriously tackled in recent years, but the bounds obtained are like

$$C_\varepsilon e^{-(1/4)(1-\varepsilon)t}, \tag{5}$$

$\varepsilon$ being an arbitrary strictly positive number with $C_\varepsilon$ going to infinity as $\varepsilon$ goes to zero. An expression of type (5) was first established in Carlen, Gabetta and Toscani (1999) under hypotheses of three different kinds on $f_0$: finiteness of all absolute moments; Sobolev regularity in the sense that $f_0$ is supposed to belong to $H^m(\mathbb{R})$ for any integer $m$; finiteness of the Linnik functional $I[f_0]$. See Linnik (1959) and Section 8 of McKean (1966) for its definition. A further progress is made in Carlen, Carvalho and Gabetta (2005), where it is shown that the above second group of hypotheses can be removed in order to obtain (5).

With a view to placing our work within the context of existing literature, the following could be to the point.

2.2. *Discussion about assumptions made in Theorem 2.1.* To start with, our moment assumption shows that the finiteness of all moments so far adopted is actually redundant. Also, the finiteness of $I[f_0]$ is not needed since in view of Lemma 2.3 in Carlen, Gabetta and Toscani (1999), $\int_\mathbb{R} e^{i\xi x} f_0(x)\,dx \leq |\xi|^{-1}\sqrt{I[f_0]}$. Hence, the tail assumption on $\varphi_0$, required in Theorem 2.1 turns out to be weaker than finiteness of $I[f_0]$.

As to the "independence" of conditions adopted in Theorem 2.1, it should be noticed that initial characteristic functions

$$\varphi_0(\xi) = \sum_{n=1}^\infty a_n \left(\frac{1}{1+\xi^2}\right)^{1/n} \qquad \left(a_n > 0 \ \forall n \text{ and } \sum_{n\geq 1} a_n = 1\right) \tag{6}$$

possess the moment property but do not meet the tail condition. Conversely, the Fourier transform of

$$f_0(x) \propto \frac{1}{1+|x|^{m+1}} \qquad \text{(for some } m \geq 1\text{)}$$

has "good" tails but $f_0$ does not possess $m$th moment.

At this stage, one could wonder whether the assumptions made in Theorem 2.1 may be weakened in some significant manner preserving, at the same time, the validity of the rate $-1/4$. Here, we explain some reasons—from (a) to (c)—why weakening of those assumptions is essentially inconsistent with that rate. (a) Theorem 4 of Carlen and Lu (2003) shows that there exist initial data $f_0$, with $m_2 < +\infty$ and $\int |x|^{2+\delta} f_0(x)\,dx = +\infty$ for every $\delta > 0$, for which $\|f(\cdot,t) - g_\sigma\|_1$ goes to zero with nonexponential rates; where finiteness of some moment of order greater than two is necessary to have bounds which decrease exponentially. (b) In any case, the following simple example seems decisive in the present discussion. Consider the class of initial densities

$$f_{0,\beta}(x) = \frac{\beta}{2|x|^{1+\beta}} \mathbb{1}_{\{|x|\geq 1\}} \qquad (x \in \mathbb{R})$$

with $3 < \beta < 4$ and Fourier transform

$$\varphi_{0,\beta}(\xi) = 1 - \frac{\beta}{2(\beta-2)}\xi^2 - \Gamma(1-\beta)\cos(\beta\pi/2)|\xi|^\beta - \beta\sum_{m\geq 2}\frac{(-1)^m \xi^{2m}}{(2m)!(2m-\beta)}.$$

Then $f_{0,\beta}$ has finite Linnik functional and entropy, finite third moment but infinite fourth moment. Moreover, $\sigma^2 = \beta(\beta-2)^{-1}$. From arguments, we will repeatedly resort to in the following Section 4, the solution $f_\beta(\cdot,t)$ of (1), with initial datum $f_{0,\beta}$, satisfies

$$\|f_\beta(\cdot,t) - g_\sigma\|_1 \geq C\exp\{-(1-2\alpha_\beta)t\} \qquad (t \geq 0)$$

for some strictly positive constant $C$ and $\alpha_\beta := \frac{1}{2\pi}\int_0^{2\pi}|\sin\theta|^\beta\,d\theta > 3/8$. This is tantamount to saying that $\|f_\beta(\cdot,t) - g_\sigma\|_1$ goes to zero exponentially, but with a rate which is slower than that provided by Theorem 2.1 under the assumption of finiteness of the fourth moment. (c) It should be noted that assumptions about finiteness of the entropy or of the Linnik functional could not compensate an eventual lack of finiteness for certain moments.

As a remark about the tail hypothesis for $\varphi_0$, one can observe that it is *slightly stronger* than the hypothesis that there is some integer $N$ such that the $N$-fold convolution of $f_0$ is bounded. This fact is of some importance since the latter assumption often recurs with a role of sufficient condition for the validity of local limit theorems in the usual Lindeberg–Lévy framework. See, for example, Sections 2 and 3 in Chapter 7 of Petrov (1975).

There is, in fact, a new problem which deserves attention. Can the upper bound we stated in Theorem 2.1 be improved? It will be shown in a forthcoming paper that the answer is in the affirmative *only if* the kurtosis coefficient of $\mu_0$ is zero, which in any case, is a rather particular condition.

We conclude the section by stressing a useful consequence of both the above discussed assumptions, which is preliminary to the proof of Theorem 2.1.

PROPOSITION 2.2. *Let $\int xf_0(x) = 0$, $\overline{m}_3 = \int |x|^3 f_0(x) < +\infty$ and assume that (3) is in force. Then there are strictly positive numbers $\lambda$ and $\alpha$ so that*

(7) $$|\varphi_0(\xi)| \leq \left(\frac{\lambda^2}{\lambda^2 + \xi^2}\right)^\alpha$$

*holds true for all real $\xi$.*

Determination of $\alpha$ and $\lambda$ are given in the course of the proof of the proposition in Appendix A.1. Notice that the right-hand side of (7) is the characteristic function of the difference $X_1 - X_2$ when $X_1$ and $X_2$ are independent and identically distributed random variables with common probability density function $\lambda^\alpha x^{\alpha-1} e^{-\lambda x} \mathbb{1}_{\{x>0\}}/\Gamma(\alpha)$.

**3. Preliminaries.** This section contains further preliminary notions and results to be used to prove Theorem 2.1. It is divided into two parts: the first explains the McKean probabilistic interpretation of *Wild's series*; the second includes estimates of the discrepancy between the normal characteristic function (the derivative of the normal characteristic function, resp.) and the characteristic function (the derivative of the characteristic function, resp.) of a sum of weighted independent and identically distributed random variables. To the best of the authors' knowledge, these estimates are new.

3.1. *Wild's sum and its probabilistic interpretation.* Let $\varphi(\cdot, t)$ be the Fourier transform of the solution $f(\cdot, t)$ of (1), that is,

$$\varphi(\xi, t) := \int_{\mathbb{R}} e^{i\xi x} f(x, t) \, dx \qquad (\xi \in \mathbb{R}, t \geq 0).$$

Then following Wild (1951), one can express $\varphi(\xi, t)$ as

$$(8) \qquad \varphi(\xi, t) := \sum_{n \geq 1} e^{-t}(1 - e^{-t})^{n-1} \hat{q}_n(\xi; \varphi_0),$$

where

$$\begin{cases} \hat{q}_1(\xi; \varphi_0) := \varphi_0(\xi) \\ \hat{q}_n(\xi; \varphi_0) = \dfrac{1}{n-1} \sum_{k=1}^{n-1} \hat{q}_k(\xi; \varphi_0) \star \hat{q}_{n-k}(\xi; \varphi_0) \qquad (n \geq 2) \end{cases}$$

and

$$g_1(\xi) \star g_2(\xi) := \frac{1}{2\pi} \int_0^{2\pi} g_1(\xi \cos \theta) \cdot g_2(\xi \sin \theta) \, d\theta.$$

For the sake of computational convenience, it is useful to recall that $\hat{q}_n$ coincides, for any $n = 2, 3, \ldots$, with the Fourier transform of the $n$-fold Wild convolution when $f_0$ is replaced by

$$(9) \qquad \tilde{f}_0(x) := \frac{f_0(x) + f_0(-x)}{2} \qquad (x \in \mathbb{R}).$$

On the one hand, the use of (9) simplifies computations required by the application of certain classical central limit arguments. On the other hand, as to the solution $\tilde{f}(\cdot, t)$ of (1) with initial condition $\tilde{f}_0$, one gets

$$(10) \qquad \int_{\mathbb{R}} |\tilde{f}(v, t) - f(v, t)| \, dv \leq 2e^{-t} \qquad (t \geq 0).$$

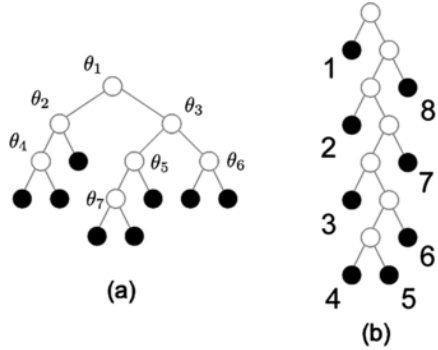

Fig. 1. *Shaded (unshaded) circles stand for leaves (nodes).*

One of the main achievements in McKean (1966) is the reinterpretation of (8) as characteristic function of a random sum. See also the recent paper [Gabetta and Regazzini (2008)]. To understand this statement, which plays a fundamental role also in the present paper, consider the product space

$$\Omega := \mathbb{N} \times \mathbb{G} \times [0, 2\pi[^\infty \times \mathbb{R}^\infty,$$

where $\mathbb{G}$ is the union $\bigcup_{n \geq 1} \mathbb{G}(n)$ and, for each $n$, $\mathbb{G}(n)$ is the set of all *McKean binary trees* with $n$ leaves [see McKean (1966, 1967)]. Each node of these trees has either zero or two "children," a "left child" and a "right child." Examples of trees are visualized in Figure 1.

Now, equip $\Omega$ with the $\sigma$-algebra

$$\mathscr{F} := \mathscr{N} \otimes \mathscr{G} \otimes \mathscr{B}([0, 2\pi[^\infty) \otimes \mathscr{B}(\mathbb{R}^\infty),$$

in which $\mathscr{N}$ and $\mathscr{G}$ are the power sets of $\mathbb{N}$ and $\mathbb{G}$, respectively, and $\mathscr{B}(S)$ denotes the Borel $\sigma$-algebra in $S$. At this stage, define

$$\nu, \tau, \theta := (\theta_n)_{n \geq 1}, \qquad \upsilon := (\upsilon_n)_{n \geq 1}$$

to be the coordinate random variables of $\Omega$.

In each tree which belongs to $\mathbb{G}(n)$, fix an order on the set of the $(n-1)$ nodes and, accordingly associate the random variable $\theta_k$ with the $k$th node. See Figure 1(a). Moreover, call $1, 2, \ldots, n$ the $n$ leaves following a left to right order. See Figure 1(b).

Define the *depth* of leaf $j$—in symbols, $\delta_j$—to be the number of generations which separate $j$ from the "root" node. With these elements at one's disposal, for each leaf $j$ of the tree $g$, one can introduce the product

$$\pi_j = \prod_{i=1}^{\delta_j} \alpha_i^{(j)},$$

where: $\alpha^{(j)}_{\delta_j}$ equals $\cos\theta_k$ if $j$ is a "left child" or $\sin\theta_k$ if $j$ is a "right child," and $\theta_k$ is the element of $\theta$ associated to the parent node of $j$; $\alpha^{(j)}_{\delta_j-1}$ equals $\cos\theta_m$ or $\sin\theta_m$ depending on the parent of $j$ is in its turn, a "left child" or a "right child," $\theta_m$ being the element of $\theta$ associated with the grandparent of $j$, and so on. For the unique tree $g_1$ in $\mathbb{G}(1)$, it is assumed that $\pi_1 \equiv 1$. From the definition of the $\pi_j$s, it is easy to deduce that

$$\sum_{j=}^{\nu} \pi_j^2 \equiv 1 \tag{11}$$

holds true whenever $\tau$ belongs to $\mathbb{G}(\nu)$. For further details of McKean's binary trees [see, e.g., McKean (1967), Carlen, Carvalho and Gabetta (2000, 2005), Gabetta and Regazzini (2006a), Bassetti, Gabetta and Regazzini (2007)].

It is easy to verify that there is one and only one probability measure $\mathsf{P}_t$ on $(\Omega, \mathscr{F})$ such that

$$\begin{aligned}&\mathsf{P}_t\{\nu = n, \tau = g, \theta \in A, \upsilon \in B\} \\ &= \begin{cases} 0, & \text{if } g \notin \mathbb{G}(n), \\ e^{-t}(1-e^{-t})^{n-1} p_n(g) u^\infty(A) \mu_0^\infty(B), & \text{if } g \in \mathbb{G}(n), \end{cases}\end{aligned} \tag{12}$$

where:

- $p_n(g)$ is the probability that a certain random walk on McKean graphs passes through the particular graph $g$ in $\mathbb{G}(n)$, $n = 1, 2, \ldots$. Such a random walk is described in McKean's papers already mentioned.
- $u^\infty$ is the probability measure on $([0, 2\pi[^\infty, \mathscr{B}([0, 2\pi[^\infty))$ that makes the $\theta_n$'s independent and identically distributed with continuous uniform law on $[0, 2\pi[$.
- $\mu_0^\infty$ is the probability measure on $(\mathbb{R}^\infty, \mathscr{B}(\mathbb{R}^\infty))$ according to which the $\upsilon_n$'s prove to be independent and identically distributed with common law $\mu_0$.

At this stage, we have all the elements needed to state the following noteworthy fact.

The solution $f(\cdot, t)$ of (1) can be viewed as a probability density function for the random variable

$$V_t = \sum_{j=1}^{\nu} \pi_j \cdot \upsilon_j \tag{13}$$

for any positive $t$.

This is tantamount to saying that $\varphi(\cdot, t)$ is the characteristic function of (13) [see McKean (1966), Gabetta and Regazzini (2008)].

It is worth mentioning that according to $\mathsf{P}_t$, the $v_n$'s turn out to be *conditionally independent* given
$$\beta := (\nu, \tau, \theta).$$
Moreover, since $(v_n)_{n\geq 1}$ and $\beta$ are independent, then in view of a standard disintegration argument, there is a conditional probability distribution $\mathsf{P}_t^*$ on $(\Omega, \mathscr{F})$ given $\beta$, according to which $V_t$ can be thought of as a sum of weighted independent random variables. See, for example, Chapter 6 in Kallenberg (2002) for existence and uniqueness results. Accordingly, $\mathsf{E}_t^*$ will denote expectation with respect to $\mathsf{P}_t^*$. It turns out that
$$\mathsf{E}_t^*[V_t] = \left(\sum_{j=1}^{\nu} \pi_j\right) \cdot \int_{\mathbb{R}} v f_0(v)\, dv$$
and via (11),
$$\mathsf{E}_t^*[V_t^2] = \int_{\mathbb{R}} v^2 f_0(v)\, dv + \left(\sum_{1\leq i\neq j\leq \nu} \pi_i \cdot \pi_j\right) \cdot \left(\int_{\mathbb{R}} v f_0(v)\, dv\right)^2.$$
Hence,
$$\mathsf{E}_t[V_t] = e^{-t} \cdot \int_{\mathbb{R}} v f_0(v)\, dv$$
and
$$\mathsf{E}_t[V_t^2] = \int_{\mathbb{R}} v^2 f_0(v)\, dv,$$
which states that the initial mean energy is preserved for any positive $t$.

The conditional independence of the $v_n$'s suggests that the study of the limiting behavior of the solution of (1), as $t$ goes to infinity, could be conducted in accordance with the ideas underlying the central limit theorem. So, the (truncated) moments of the summands $\pi_j v_j$ play a fundamental role with respect to both conditions for convergence and precise estimation of the rates of convergence. Since these moments depend essentially on sums of powers of $\pi_j$'s, the following identities, proved in Gabetta and Regazzini (2006a), deserve attention:
$$\mathsf{E}_t\left[\sum_{j=1}^{\nu} |\pi_j|^m \Big| \nu\right] = \frac{\Gamma(2\alpha_m + \nu - 1)}{\Gamma(2\alpha_m)\Gamma(\nu)}$$
and
$$(14) \qquad \mathsf{E}_t\left[\sum_{j=1}^{\nu} |\pi_j|^m\right] = e^{-(1-2\alpha_m)t},$$
where
$$(15) \qquad \alpha_m := \frac{1}{2\pi} \int_0^{2\pi} |\sin\theta|^m\, d\theta.$$

3.2. *Measures of discrepancy between characteristic functions of sums and normal characteristic function.* To pave the way for the application of some classical probabilistic results to the study of the rate of convergence of the solution of (1), it is essential to fit a few well-known statements, on characteristic function of sums of independent random variables, to the setup determined by the peculiarity of the sum $V_t$ in (13). In the remainder of this subsection, $X_1, X_2, \ldots, X_n$ will stand for independent and identically distributed real-valued random variables on $(\tilde{\Omega}, \tilde{\mathscr{F}}, \tilde{\mathsf{P}})$, with common non-degenerate distribution $\tilde{\mu}_0$. Moreover, assume that $\tilde{\mu}_0$ is symmetric and has moment of fourth power. Denote the $k$th moment and the absolute $k$th moment of $\tilde{\mu}_0$ by $m_k$ and $\overline{m}_k$, respectively. Notice that $m_1 = m_3 = 0$. Set $\tilde{\varphi}_0(\xi) := \int_{\mathbb{R}} e^{i\xi x} d\tilde{\mu}_0(x)$, and let $c_1, \ldots, c_n$ be real constants such that

$$\sum_{j=1}^{n} c_j^2 = 1.$$

Then define $V_n$ to be the sum of $Y_1, \ldots, Y_n$, with

$$Y_j := \frac{1}{\sqrt{m_2}} c_j X_j \qquad (j = 1, \ldots, n)$$

and let $\tilde{\varphi}_n$ be the characteristic function of $V_n$. Finally, set

$$\Gamma_n^4 := \frac{m_4}{m_2^2} \sum_{j=1}^{n} c_j^4.$$

LEMMA 3.1. *The following inequalities*

(16) $$\left( \int_{-A}^{A} |\tilde{\varphi}_n(\xi) - e^{-\xi^2/2}|^2 \, d\xi \right)^{1/2} \leq 0.62 \sqrt{\Gamma(7/2)} \Gamma_n^4,$$

(17) $$\left( \int_{-A}^{A} \left| \frac{d}{d\xi} \{\tilde{\varphi}_n(\xi) - e^{-\xi^2/2}\} \right|^2 d\xi \right)^{1/2} \leq 3.8 \sqrt{\Gamma(7/2)} \Gamma_n^4$$

*are valid for $A = a/\Gamma_n$, $a$ being any positive number not greater than $1/2$.*

Clearly, this proposition is reminiscent of the Berry–Esseen inequality. In fact, the proof of Lemma 3.1, deferred to Section A.2 of the Appendix, is based on the exposition of that inequality to be found in Chow and Teicher (1997).

**4. Proof of Theorem 2.1.** From now on, $\mathsf{F}^*(x)$ will stand for conditional probability $\mathsf{P}_t^* \{V_t \leq x\}$. Moreover, integral over a measurable subset $S$ of $\Omega$ will be often denoted by $\mathsf{E}[\cdot; S]$. The notation $m_k$ and $\overline{m}_k$ for $\int v^k f_0(v)$ and $\int |v|^k f_0(v)$, respectively, will also be as a norm resorted to.

To prove the theorem, first recall definition (9) and inequality (10) to guarantee, through a plain application of the triangle inequality, that it is enough to verify the validity of inequality

$$\int_{\mathbb{R}} |\tilde{f}(v,t) - g_\sigma(v)| \, dv \leq C_* e^{-(1/4)t} \tag{18}$$

for some suitable $C_*$. Working with $\tilde{f}(v,t)$, instead of $f(v,t)$, has the advantage that the existing odd moments of $\tilde{f}(\cdot,t)$ vanish. In view of this, for the sake of notational convenience, it will be assumed that condition

(H)    $f_0$ and, consequently, $f(\cdot,t)$ are even functions

is in force. Clearly,

$$\|f(v,t) - g_\sigma(v)\|_1 \leq \mathsf{E}_t\left[\left\|\frac{d}{dv}\mathsf{F}^*(v) - g_\sigma(v)\right\|_1\right]$$

$$= \mathsf{E}_t\left[\left\|\frac{d}{dv}\mathsf{F}^*(\sigma v) - g_1(v)\right\|_1\right]. \tag{19}$$

Now, notice that Proposition 2.2 is applicable to $f_0$, thanks to the hypotheses of Theorem 2.1 and (H). Therefore, once $\lambda$ and $\alpha$ have been fixed as in the above mentioned proposition, one can define $U \subset \Omega$ as

$$U := \{\nu \leq \overline{n}\} \cup \left\{\prod_{j=1}^{\nu} \pi_j = 0\right\} \cup \left\{\sum_{j=1}^{\nu} \pi_j^4 \geq \overline{\delta}\right\} \tag{20}$$

$\overline{n}$ being equal to $\lceil 9/2\alpha \rceil$ and $\overline{\delta} = (2^{\overline{n}} \overline{n}!)^{-1}$. Next, check that $U$ belongs to $\mathscr{F}$ and rewrite the last term in (19) as

$$\mathsf{E}_t\left[\left\|\frac{d}{dv}\mathsf{F}^*(\sigma v) - g_1(v)\right\|_1\right]$$

$$= \mathsf{E}_t\left[\left\|\frac{d}{dv}\mathsf{F}^*(\sigma v) - g_1(v)\right\|_1; U\right] + \mathsf{E}_t\left[\left\|\frac{d}{dv}\mathsf{F}^*(\sigma v) - g_1(v)\right\|_1; U^C\right]. \tag{21}$$

In order to obtain estimates of these terms, observe that

$$\mathsf{P}_t(U) \leq (\overline{n} + 2^{\overline{n}}\overline{n}!)e^{-(1/4)t} \tag{22}$$

is valid for all nonnegative $t$. Indeed,

$$\mathsf{P}_t\{\nu \leq \overline{n}\} = \sum_{k=1}^{\overline{n}} e^{-t}(1 - e^{-t})^{k-1} \leq \overline{n}e^{-t} \leq \overline{n}e^{-(1/4)t}.$$

Moreover, as a consequence of the fact that $\mathsf{P}_t\{\prod_{j=1}^{\nu} \pi_j = 0 \mid \nu = m, \tau = g\}$ is 0 for all $(m,g)$ in $\mathbb{N} \times \mathbb{G}$, one gets

$$\mathsf{P}_t\left\{\prod_{j=1}^{\nu} \pi_j = 0\right\} = 0.$$

Finally, one can put $m = 4$ in (15) to have $\alpha_4 = 3/8$, and combine the Markov inequality with (14) to obtain

$$\mathsf{P}_t\left\{\sum_{j=1}^{\nu}\pi_j^4 \geq \overline{\delta}\right\} \leq \frac{1}{\overline{\delta}}e^{-(1/4)t} = 2^{\overline{n}}\overline{n}!e^{-(1/4)t}$$

from which (22) follows.

On the one hand, from $\|\frac{d}{dv}\mathsf{F}^*(\sigma v) - g_1(v)\|_1 \leq 2$ one deduces the following upper bound:

$$(23) \quad \mathsf{E}_t\left[\left\|\frac{d}{dv}\mathsf{F}^*(\sigma v) - g_1(v)\right\|_1; U\right] \leq 2\mathsf{P}_t(U) \leq 2(\overline{n} + 2^{\overline{n}}\overline{n}!)e^{-(1/4)t}.$$

On the other hand, since $U^c$ can be viewed as a "good" set, it is to be expected that $\|\frac{d}{dv}\mathsf{F}^*(\sigma v) - g_1(v)\|_1$ is very small on $U^c$. To check this fact, one can resort to Proposition 4.1.

PROPOSITION 4.1 [Beurling (1939)]. *Let $f$ be in $L^1(\mathbb{R})$, that is,*

$$\int_{\mathbb{R}} |f(x)|\, dx < +\infty$$

*with Fourier transform $\varphi$ in the Sobolev space $H^1(\mathbb{R})$, that is:*

- $\int_{\mathbb{R}} |\varphi(\xi)|^2\, d\xi < +\infty$;
- *there exists an essentially unique function $\psi$ such that $\int_{\mathbb{R}} |\psi(\xi)|^2\, d\xi < +\infty$ and*

$$\int_{\mathbb{R}} \varphi(\xi)\phi'(\xi)\, d\xi = -\int_{\mathbb{R}} \psi(\xi)\phi(\xi)\, d\xi \qquad (\phi \in C_c^\infty(\mathbb{R})).$$

*Then*

$$(24) \quad \sqrt{2}\|f\|_1 \leq \|\varphi\|_{H^1(\mathbb{R})} := \left(\int_{\mathbb{R}} |\varphi(\xi)|^2\, d\xi + \int_{\mathbb{R}} |\psi(\xi)|^2\, d\xi\right)^{1/2}.$$

For the sake of completeness, a proof of this proposition—briefly sketched in Beurling (1939)—will be given in Section A.1 in the Appendix.

Under the assumptions of Theorem 2.1, the restriction to $U^c$ of the conditional characteristic function $\xi \mapsto \varphi^*(\xi) := \int_{\mathbb{R}} e^{i\xi x}\, d\mathsf{F}^*(x)$ belongs to $H^1(\mathbb{R})$—see Remarks A.1 and A.2 in the Appendix—and this, in view of Proposition 4.1, yields

$$\mathsf{E}_t\left[\left\|\frac{d}{dv}\mathsf{F}^*(\sigma v) - g_1(v)\right\|_1; U^c\right] \leq \frac{1}{\sqrt{2}}\mathsf{E}_t\left[\left\{\int_{\mathbb{R}} |\Delta|^2\, d\xi + \int_{\mathbb{R}} |\Delta'|^2\, d\xi\right\}^{1/2}; U^c\right],$$

where $\Delta := \varphi^*(\frac{\xi}{\sigma}) - e^{-\xi^2/2}$ and $\Delta' := \frac{d}{d\xi}\Delta$. Now,

$$\mathsf{E}_t\left[\left\{\int_{\mathbb{R}}|\Delta|^2\,d\xi + \int_{\mathbb{R}}|\Delta'|^2\,d\xi\right\}^{1/2}; U^c\right]$$

(25)
$$\leq \mathsf{E}_t\left[\left(\int_{\{|\xi|\leq A\}}|\Delta|^2\,d\xi\right)^{1/2}; U^c\right] + \mathsf{E}_t\left[\left(\int_{\{|\xi|\geq A\}}|\Delta|^2\,d\xi\right)^{1/2}; U^c\right]$$

$$+ \mathsf{E}_t\left[\left(\int_{\{|\xi|\leq A\}}|\Delta'|^2\,d\xi\right)^{1/2}; U^c\right]$$

$$+ \mathsf{E}_t\left[\left(\int_{\{|\xi|\geq A\}}|\Delta'|^2\,d\xi\right)^{1/2}; U^c\right]$$

with

$$A = A(\beta) := \frac{\sigma^4}{2m_4(\sum_{j=1}^{\nu}\pi_j^4)^{1/4}}.$$

Thanks to the Lyapunov inequality one has

$$A \leq \frac{\sigma}{2m_4^{1/4}(\sum_{j=1}^{\nu}\pi_j^4)^{1/4}} = \frac{1}{2\Gamma_\nu(\pi_1,\ldots,\pi_\nu)};$$

therefore, one can apply Lemma 3.1 to obtain

$$\left(\int_{\{|\xi|\leq A\}}|\Delta|^2\,d\xi\right)^{1/2} \leq 0.62\sqrt{\Gamma(7/2)}\frac{m_4}{\sigma^4}\cdot\sum_{j=1}^{\nu}\pi_j^4$$

and

$$\left(\int_{\{|\xi|\leq A\}}|\Delta'|^2\,d\xi\right)^{1/2} \leq 3.8\sqrt{\Gamma(7/2)}\frac{m_4}{\sigma^4}\cdot\sum_{j=1}^{\nu}\pi_j^4.$$

At this stage, take expectation $\mathsf{E}_t$ and recall (14) with $m = 4$ ($\alpha_4 = 3/8$) to obtain both

(26) $$\mathsf{E}_t\left[\left(\int_{\{|\xi|\leq A\}}|\Delta|^2\,d\xi\right)^{1/2}\right] \leq 0.62\sqrt{\Gamma(7/2)}\frac{m_4}{\sigma^4}\cdot e^{-(1/4)t}$$

and

(27) $$\mathsf{E}_t\left[\left(\int_{\{|\xi|\leq A\}}|\Delta'|^2\,d\xi\right)^{1/2}\right] \leq 3.8\sqrt{\Gamma(7/2)}\frac{m_4}{\sigma^4}\cdot e^{-(1/4)t}.$$

After determining upper bounds for integrals of the type of $\int_{\{|\xi|\leq A\}}$, one gets down to examining the remaining summands in (25). The Minkowski

inequality yields

$$\left(\int_{\{|\xi|\geq A\}}|\Delta|^2\,d\xi\right)^{1/2}$$
$$\leq \left(\int_{\{|\xi|\geq A\}}|\varphi^*(\xi/\sigma)|^2\,d\xi\right)^{1/2} + \left(\int_{\{|\xi|\geq A\}}|e^{-\xi^2/2}|^2\,d\xi\right)^{1/2}$$

and

$$\left(\int_{\{|\xi|\geq A\}}|\Delta'|^2\,d\xi\right)^{1/2}$$
$$\leq \left(\int_{\{|\xi|\geq A\}}\left|\frac{d}{d\xi}\varphi^*(\xi/\sigma)\right|^2\,d\xi\right)^{1/2} + \left(\int_{\{|\xi|\geq A\}}|\xi e^{-\xi^2/2}|^2\,d\xi\right)^{1/2}.$$

Combining a well-known inequality, proved, for example, in Lemma 2 of VII.1 in Feller (1968), with $\max_{x\geq 0} x^k e^{-\alpha x^2} = [k/(2e\alpha)]^{k/2}$, one obtains

$$\left(\int_{\{|\xi|\geq A\}} e^{-\xi^2}\,d\xi\right)^{1/2} \leq 16 e^{-2}\left(\frac{2m_4}{\sigma^4}\right)^{9/2}\sum_{j=1}^{\nu}\pi_j^4$$

and

$$\left(\int_{\{|\xi|\geq A\}} \xi^2 e^{-\xi^2}\,d\xi\right)^{1/2} \leq c\left(\frac{2m_4}{\sigma^4}\right)^{9/2}\sum_{j=1}^{\nu}\pi_j^4,$$

where $c$ can be chosen equal to $(5/e)^{5/2} + 8\sqrt{2}e^{-2}$. Then from the last inequalities together with (14), it follows,

$$(28) \qquad \mathsf{E}_t\left(\int_{\{|\xi|\geq A\}} e^{-\xi^2}\,d\xi\right)^{1/2} \leq 16 e^{-2}\left(\frac{2m_4}{\sigma^4}\right)^{9/2} e^{-(1/4)t}$$

and

$$(29) \qquad \mathsf{E}_t\left(\int_{\{|\xi|\geq A\}} \xi^2 e^{-\xi^2}\,d\xi\right)^{1/2} \leq c\left(\frac{2m_4}{\sigma^4}\right)^{9/2} e^{-(1/4)t}.$$

From Remark A.2 in the Appendix,

$$\left[\left(\int_{\{|\xi|\geq A\}}|\varphi^*(\xi/\sigma)|^2\,d\xi\right)^{1/2} + \left(\int_{\{|\xi|\geq A\}}\left|\frac{d}{d\xi}\varphi^*(\xi/\sigma)\right|^2\,d\xi\right)^{1/2}\right]\cdot \mathbb{1}_{U^C}$$
$$= \sqrt{2}\left[\left(\int_A^{+\infty}|\varphi^*(\xi/\sigma)|^2\,d\xi\right)^{1/2}\right.$$
$$(30)$$
$$\left. + \left(\int_A^{+\infty}\left|\frac{d}{d\xi}\varphi^*(\xi/\sigma)\right|^2\,d\xi\right)^{1/2}\right]\cdot \mathbb{1}_{U^C}$$
$$\leq 2\sqrt{2}\left(\int_A^{+\infty}|\varphi^*(\xi/\sigma)|\,d\xi\right)^{1/2}\cdot \mathbb{1}_{U^C} + \sqrt{2|\varphi^*(A/\sigma)|}.$$

Now one can resort to the classical argument used to prove the Berry–Esseen inequality—see, for example, Lemma 12 in Chapter VI of Petrov (1975); the applicability of this result (with $b = 1/2$) is guaranteed by the inequality

$$A \leq \frac{\sigma^3}{2\overline{m}_3 \sum_{j=1}^{\nu} |\pi_j|^3}.$$

Thus,

$$\sqrt{2|\varphi^*(A/\sigma)|} \leq \sqrt{2}e^{-(1/12)A^2} \leq \sqrt{2}(24/e)^2 A^{-4} = \sqrt{2}(24/e)^2 \left(\frac{2m_4}{\sigma^4}\right)^4 \sum_{j=1}^{\nu} \pi_j^4$$

and, therefore,

$$(31) \qquad \mathsf{E}_t\sqrt{2|\varphi^*(A/\sigma)|} \leq \sqrt{2}(24/e)^2 \left(\frac{2m_4}{\sigma^4}\right)^4 e^{-(1/4)t}.$$

There is one more thing to analyze, that is,

$$(32) \quad \left(\int_A^{+\infty} |\varphi^*(\xi/\sigma)|\, d\xi\right)^{1/2} \cdot \mathbb{1}_{U^C} = \left(\int_A^{+\infty} \prod_{j=1}^{\nu} \left|\varphi_0\left(\frac{\pi_j \xi}{\sigma}\right)\right| d\xi\right)^{1/2} \cdot \mathbb{1}_{U^C}.$$

In view of Proposition 2.2, the right-hand side turns out to be not greater than

$$(33) \quad \begin{aligned} & \left[\int_A^{+\infty} \prod_{j=1}^{\nu} \left(\frac{\lambda^2}{\lambda^2 + ((\pi_j \xi)/\sigma)^2}\right)^{\alpha} d\xi\right]^{1/2} \cdot \mathbb{1}_{U^C} \\ &= \left[\lambda \sigma \int_{A/\lambda\sigma}^{+\infty} \left(\frac{1}{\prod_{j=1}^{\nu}(1 + \pi_j^2 \eta^2)}\right)^{\alpha} d\eta\right]^{1/2} \cdot \mathbb{1}_{U^C}.\end{aligned}$$

Moreover, inequality

$$(34) \qquad \prod_{j=1}^{\nu} (1 + \pi_j^2 \eta^2) \geq \overline{\varepsilon} \eta^{2\overline{n}}$$

is valid on $U^c$, with $\overline{\varepsilon} := 1/\overline{n}! - 2^{\overline{n}-1} \cdot \overline{\delta} = 1/2\overline{n}!$. See Section A.3 in the Appendix. Then

$$(35) \quad \begin{aligned} & \left[\lambda \sigma \int_{A/\lambda\sigma}^{+\infty} \left(\frac{1}{\prod_{j=1}^{\nu}(1+\pi_j^2\eta^2)}\right)^{\alpha} d\eta\right]^{1/2} \cdot \mathbb{1}_{U^C} \\ & \leq \left[\lambda \sigma \int_{A/\lambda\sigma}^{+\infty} \left(\frac{1}{\overline{\varepsilon}\eta^{2\overline{n}}}\right)^{\alpha} d\eta\right]^{1/2} = \tilde{C}\left(\sum_{j=1}^{\nu} \pi_j^4\right)^{(2\overline{n}\alpha-1)/8}.\end{aligned}$$

The definition of $\overline{n}$ in (20) yields equality $(2\alpha\overline{n} - 1)/8 = 1$. Moreover,

$$
\begin{aligned}
(36) \quad \tilde{C} &:= \frac{(\lambda\sigma)^{9/2}}{2\sqrt{2}\overline{\varepsilon}^{\alpha/2}}\left(\frac{2m_4}{\sigma^4}\right)^4 = 4\sqrt{2}\left(\frac{m_4^4}{\sigma^{23/2}}\right)(2\overline{n}!)^{9/(4\overline{n})}\lambda^{9/2} \\
&\leq 4\sqrt{2}\left(\frac{m_4^4}{\sigma^{23/2}}\right)(2\overline{n}!)^{9/(4\overline{n})}2^{5/4}\left[\left(\frac{3}{2\sigma^2}\right)^{9/4} + \left(\frac{2}{1-M}\right)^{9/4}(L_p)^{9/2p}\right]
\end{aligned}
$$

with

$$
M = \exp\left\{-\frac{3\pi^2}{64(3 + (L_p)^{4/p})^2}\left(\frac{\sqrt{2}\sigma}{8\lceil 2/p\rceil\sigma^3 + 40\pi\sqrt{\lceil 2/p\rceil}m_4}\right)^2\right\}.
$$

Hence, (35) and (14) yield

$$
(37) \quad \mathsf{E}_t\left[\left(\int_A^{+\infty}|\varphi^*(\xi/\sigma)|\,d\xi\right)^{1/2}\cdot \mathbb{1}_{U^C}\right] \leq \tilde{C}e^{-(1/4)t}
$$

and, to complete the proof of the theorem, it is enough to combine (23), (26)–(29), (31) and (37).

## APPENDIX

For the sake of completeness, we gather here certain propositions required in the text, as well as a few proofs skipped in the previous sections. The subject is split into three subsections. The first subsection includes proofs for Propositions 2.2 and 4.1. The second contains adaptations to the specific setting of the present paper of well-known results pertaining to asymptotic expansions in the central limit theorem. Finally, the last one completes the proof of Theorem 2.1 provided in Section 2.

**A.1.** A preliminary lemma is needed.

LEMMA A.1. *For any probability measure $\nu$ on $(\mathbb{R}, \mathscr{B}(\mathbb{R}))$ with the Fourier–Stieltjes transform $\psi(\xi) := \int e^{i\xi x}\,d\nu(x)$, such that:*

- $\int x\,d\nu(x) = 0$,
- $\int x^2\,d\nu(x) =: \zeta^2 < +\infty$,
- $\psi(\xi) = o(|\xi|^{-4}), (|\xi| \to +\infty)$,

*one has*

$$
(38) \quad |\psi(\xi)| \leq \exp\left\{-\frac{3\pi^2}{64(3 + L)^2}\left(\frac{\xi}{2\sqrt{2}\zeta|\xi| + \pi}\right)^2\right\} \quad (\xi \in \mathbb{R}),
$$

*where $L := \sup_{\xi \in \mathbb{R}}|\xi^4\psi(\xi)|$.*

PROOF. From the tail assumption on $\psi$ and the definition of $L$, it is plain that $|\psi(\xi)| \leq \min\{1; L\xi^{-4}\}$. Thus,

$$\int_{\mathbb{R}} |\psi(\xi)|^2 \, d\xi \leq \int_{\mathbb{R}} |\psi(\xi)| \, d\xi \leq 2\left(1 + \int_1^{+\infty} \frac{L}{\xi^4} \, d\xi\right) = \frac{2}{3}(3+L).$$

As $\psi$ has finite integral, $\nu$ turns out to be absolutely continuous with density $f$ given by

$$f(x) = \frac{1}{2\pi} \int_{\mathbb{R}} \psi(\xi) e^{-i\xi x} \, d\xi.$$

Furthermore, $f$ belongs to $L^2(\mathbb{R})$, as an immediate consequence of the Plancherel identity and the fact that $\psi$ itself belongs to $L^2(\mathbb{R})$. Then

$$\int_{\mathbb{R}} f^2(x) \, dx = \frac{1}{2\pi} \int_{\mathbb{R}} \psi^2(\xi) \, d\xi \leq \frac{1}{3\pi}(3+L).$$

Now, Young's inequality for convolutions [see, e.g., Proposition 8.8 in Folland (1999)] entails

$$\sup_{x \in \mathbb{R}} |f * f(x)| \leq \int_{\mathbb{R}} f^2(x) \, dx \leq \frac{1}{3\pi}(3+L).$$

At this stage, invoke formula (6.73) on page 172 of Saulis and Statulevičius (1991) to obtain

$$|\psi(\xi)|^2 \leq \exp\left\{-\frac{9\pi^2}{96(3+L)^2}\left(\frac{\xi}{2\sqrt{2}\zeta|\xi|+\pi}\right)^2\right\} \qquad (\xi \in \mathbb{R}). \qquad \square$$

PROOF OF PROPOSITION 2.2. Define $\psi_0(\xi) := |\varphi_0(\xi)|^{2k}$ with $k = \lceil 2/p \rceil$ and notice that $\psi_0$ is a real-valued characteristic function such that, if $Y$ is any random variable with such a characteristic function, then

$$\mathsf{E}Y = 0, \qquad \mathsf{E}Y^2 = 2k\sigma^2, \qquad \mathsf{E}|Y|^3 \leq 10\sqrt{2}k^{3/2}\overline{m}_3.$$

As to the last inequality, see Petrov (1975), page 60. According to the argument used to prove the Berry–Esseen inequality [see, e.g., Petrov (1975), page 110], one gets

$$|\psi_0(\xi)|^2 \leq 1 - \frac{2}{3} 2k\sigma^2 \xi^2 = 1 - \frac{4k}{3}\sigma^2 \xi^2$$

provided that $|\xi| \leq \frac{\sqrt{2}\sigma^2}{40\sqrt{k}\overline{m}_3}$. Hence, under this very same restriction,

$$|\psi_0(\xi)| \leq \sqrt{1 - \frac{4k}{3}\sigma^2 \xi^2} \leq 1 - \frac{2k}{3}\sigma^2 \xi^2$$

and

$$1 - \frac{2k}{3}\sigma^2 \xi^2 \leq \frac{\lambda^2}{\lambda^2 + \xi^2},$$

whenever $\lambda^2 \geq \frac{3}{2k\sigma^2}$. After examining what happens to internal values of $\xi$, one can observe that

$$|\psi_0(\xi)| = o(|\xi|^{-4}) \qquad (|\xi| \to +\infty).$$

Then

$$|\psi_0(\xi)| \leq \frac{L}{\xi^4} \leq \frac{\lambda^2}{\lambda^2 + \xi^2},$$

where the last inequality holds for any $\lambda \neq 0$ and

$$\xi^2 \geq \frac{L + \sqrt{L^2 + 4L\lambda^4}}{2\lambda^2} \geq \sqrt{L}$$

with $L := \sup_{\xi \in \mathbb{R}} |\xi^4 \psi_0(\xi)|$. At this stage, the proof would be complete if $L^{1/4}$ were less than $\frac{\sqrt{2}\sigma^2}{40\sqrt{k\overline{m}_3}}$. It remains to analyze the case in which $L^{1/4} \geq \frac{\sqrt{2}\sigma^2}{40\sqrt{k\overline{m}_3}}$. If $\lambda^2 \geq (2\sqrt{L})/3$, then

$$\sqrt{L} \leq \frac{L + \sqrt{L^2 + 4L\lambda^4}}{2\lambda^2} \leq 2\sqrt{L}.$$

In view of the tail condition, it turns out that the maximum value $M$ of $|\psi_0(\xi)|$ is smaller that 1 on $I := [\frac{\sqrt{2}\sigma^2}{40\sqrt{k\overline{m}_3}}, \sqrt{2}L^{1/4}]$. Thus, in order that

$$|\psi_0(\xi)| \leq \frac{\lambda^2}{\lambda^2 + \xi^2} \tag{39}$$

be valid on $I$ it is enough to choose $\lambda^2 \geq \frac{2M\sqrt{L}}{1-M}$. Now, if

$$\lambda^2 = \max\left\{\frac{3}{2k\sigma^2}, \frac{2}{3}\sqrt{L}, \frac{2M\sqrt{L}}{1-M}\right\} \leq \frac{3}{2\sigma^2} + \frac{2}{1-M}\sqrt{L},$$

then (39) is valid for every $\xi$, and one can take $\alpha = \frac{1}{2k}$ to get (7). □

PROOF OF PROPOSITION 4.1 (Beurling). From the Cauchy–Schwarz–Buniakowsky inequality,

$$\int_\mathbb{R} |f(x)|\, dx = \int_\mathbb{R} |f(x)| \frac{\sqrt{1+x^2}}{\sqrt{1+x^2}}\, dx$$

$$\leq \left\{\int_\mathbb{R} |f(x)|^2 (1+x^2)\, dx\right\}^{1/2} \left\{\int_\mathbb{R} \frac{dx}{1+x^2}\right\}^{1/2}$$

$$= \sqrt{\pi} \left\{\int_\mathbb{R} |f(x)|^2\, dx + \int_\mathbb{R} |xf(x)|^2\, dx\right\}^{1/2}.$$

At this stage, an application of the Plancherel theorem yields

$$\sqrt{\pi}\left\{\int_{\mathbb{R}}|f(x)|^2\,dx+\int_{\mathbb{R}}|xf(x)|^2\,dx\right\}^{1/2}$$
$$=\sqrt{\pi}\left\{\frac{1}{2\pi}\int_{\mathbb{R}}|\varphi(\xi)|^2\,d\xi+\frac{1}{2\pi}\int_{\mathbb{R}}|\psi(\xi)|^2\,d\xi\right\}^{1/2}\leq\frac{1}{\sqrt{2}}\|\varphi\|_{H^1(\mathbb{R})}.\quad\square$$

Hence, for any pair of probability density functions, $f_1$ and $f_2$, having finite expectations and Fourier transforms $\varphi_1$ and $\varphi_2$ in $H^1(\mathbb{R})$, one gets

$$\int_{\mathbb{R}}|f_1(x)-f_2(x)|\,dx\leq\left\{\int_{\mathbb{R}}|\varphi_1(\xi)-\varphi_2(\xi)|^2\,d\xi+\int_{\mathbb{R}}|\varphi_1'(\xi)-\varphi_2'(\xi)|^2\,d\xi\right\}^{1/2}.$$

This inequality is explicitly mentioned in Zolotarev (1997).

**A.2.** This section aims at providing proofs for inequalities (16)–(17). To this purpose, a further preliminary result is needed anyway.

LEMMA A.2. *For $|\xi|\leq 1/(2\Gamma_n)$, one has*

$$|\tilde{\varphi}_n(\xi)-e^{-\xi^2/2}|\leq c_1\Gamma_n^4\xi^4 e^{-\xi^2/2} \tag{40}$$

*and*

$$\left|\frac{d}{d\xi}[\tilde{\varphi}_n(\xi)-e^{-\xi^2/2}]\right|\leq c_2\Gamma_n^4(1+\xi^2)|\xi|^3 e^{-\xi^2/2}. \tag{41}$$

*One can take $c_1=0.33$ and $c_2=0.76$.*

PROOF. Setting $\tilde{\varphi}_{n,j}$ for the characteristic function of $Y_{n,j}$ ($j=1,2,\ldots,n$), one has

$$\tilde{\varphi}_{n,j}(\xi)=\tilde{\varphi}_0\left(\frac{c_j}{\sqrt{m_2}}\xi\right)=1-\frac{1}{2}c_j^2\xi^2+R_4^{(j)}(\xi)$$

with

$$|R_4^{(j)}(\xi)|\leq\frac{1}{24}\mathsf{E}(Y_j^4)\xi^4=\frac{1}{24}\frac{m_4}{m_2^2}c_j^4\xi^4$$

and

$$|c_j\xi|\leq\left(\frac{m_4}{m_2^2}\right)^{1/4}|c_j\xi|\leq|\Gamma_n\xi|\leq 1/2.$$

Moreover,

$$\left|-\frac{1}{2}c_j^2\xi^2+R_4^{(j)}(\xi)\right|\leq\frac{1}{8}+\frac{1}{24}\frac{m_4}{m_2^2}c_j^4\frac{1}{2^4\Gamma_n^4}\leq\frac{49}{384},$$

which entails $335/384 \leq \tilde{\varphi}_{n,j}(\xi) \leq 1$, and

$$\begin{aligned}\log \tilde{\varphi}_{n,j}(\xi) &= \log\{1 + (-\tfrac{1}{2}c_j^2\xi^2 + R_4^{(j)}(\xi))\} \\ &= -\tfrac{1}{2}c_j^2\xi^2 + R_4^{(j)}(\xi) + (-\tfrac{1}{2}c_j^2\xi^2 + R_4^{(j)}(\xi))^2 \rho_j(\xi),\end{aligned}$$

where $\rho_j(\xi) := \rho(-\tfrac{1}{2}c_j^2\xi^2 + R_4^{(j)}(\xi))$ and

$$\rho(z) := \begin{cases} z^{-2}[\log(1+z) - z], & z > -1, z \neq 0, \\ -1/2, & z = 0. \end{cases}$$

Hence,

$$\begin{aligned}\log \tilde{\varphi}_n(\xi) &= \log\left\{\prod_{j=1}^n \tilde{\varphi}_{n,j}(\xi)\right\} = \sum_{j=1}^n \log \tilde{\varphi}_{n,j}(\xi) \\ &= -\tfrac{1}{2}\xi^2 + \sum_{j=1}^n R_4^{(j)}(\xi) + \sum_{j=1}^n (-\tfrac{1}{2}c_j^2\xi^2 + R_4^{(j)}(\xi))^2 \rho_j(\xi).\end{aligned}$$

Set

$$R_{1,n}(\xi) := \sum_{j=1}^n R_4^{(j)}(\xi),$$

$$R_{2,n}(\xi) := \sum_{j=1}^n (-\tfrac{1}{2}c_j^2\xi^2 + R_4^{(j)}(\xi))^2 \rho_j(\xi),$$

$$R_n(\xi) := R_{1,n}(\xi) + R_{2,n}(\xi).$$

Then letting $\rho^* := -\rho(-49/384) < 0,55$, one gets

$$|R_{1,n}(\xi)| \leq \sum_{j=1}^n |R_4^{(j)}(\xi)| \leq \tfrac{1}{24}\Gamma_n^4 \xi^4,$$

$$|R_{1,n}(\xi)| \leq 0.002$$

and

$$|R_{2,n}(\xi)| \leq 2\rho^* \sum_{j=1}^n [\tfrac{1}{4}c_j^4\xi^4 + (R_4^{(j)}(\xi))^2] \leq 0.29\Gamma_n^4 \xi^4,$$

$$|R_{2,n}(\xi)| \leq 0.02,$$

hence

(42)     $|R_n(\xi)| \leq 0.33\Gamma_n^4 \xi^4,$

(43)     $|R_n(\xi)| \leq 0.022.$

Now, from the well-known Bernstein form for the error term in Taylor's formula [see, e.g., Theorem 9.29 in Apostol (1974)], one obtains

$$\left|\frac{d}{d\xi}R_4^{(j)}(\xi)\right| \leq \frac{1}{6}\frac{m_4}{m_2^2}c_j^4|\xi|^3$$

and

$$\left|\frac{d}{d\xi}R_{1,n}(\xi)\right| \leq \sum_{j=1}^n \left|\frac{d}{d\xi}R_4^{(j)}(\xi)\right| \leq \frac{1}{6}\Gamma_n^4|\xi|^3.$$

Next,

$$\left|\frac{d}{d\xi}R_{2,n}(\xi)\right|$$

$$\leq \left(\frac{49}{384}\rho_1^* + 2\rho^*\right) \cdot \sum_{j=1}^n \left|\left(-\frac{1}{2}c_j^2\xi^2 + R_4^{(j)}(\xi)\right) \cdot \left(-c_j^2\xi + \frac{d}{d\xi}R_4^{(j)}(\xi)\right)\right|$$

with $\rho_1^* := \frac{d}{dz}\rho(z)|_{z=-49/384}$. Hence,

$$\left|\frac{d}{d\xi}R_{2,n}(\xi)\right| \leq 0.6\Gamma_n^4(1+\xi^2)|\xi|^3$$

and

(44) $$\left|\frac{d}{d\xi}R_n(\xi)\right| \leq 0.8\Gamma_n^4(1+\xi^2)|\xi|^3.$$

Finally, combination of (42), (43) and (44) with the elementary inequality $|e^z - 1| \leq |z|e^{|z|}$ gives

$$|\tilde{\varphi}_n(\xi) - e^{-\xi^2/2}| \leq 0.33\Gamma_n^4\xi^4 e^{-\xi^2/2}$$

and

$$\left|\frac{d}{d\xi}[\tilde{\varphi}_n(\xi) - e^{-\xi^2/2}]\right| \leq e^{|R_n(\xi)|}\left[|\xi R_n(\xi)| + \left|\frac{d}{d\xi}R_n(\xi)\right|\right]e^{-\xi^2/2}$$

$$\leq 0.76\Gamma_n^4(1+\xi^2)|\xi|^3 e^{-\xi^2/2}. \qquad \square$$

PROOF OF INEQUALITIES (16)–(17). Take $A = 1/2\Gamma_n$. From Lemma A.2,

$$\left(\int_{-A}^A |\tilde{\varphi}_n(\xi) - e^{-\xi^2/2}|^2\, d\xi\right)^{1/2} \leq 0.33\Gamma_n^4 \left(\int_{-A}^A \xi^8 e^{-\xi^2}\, d\xi\right)^{1/2}$$

$$\leq 0.62\sqrt{\Gamma(7/2)}\,\Gamma_n^4$$

and

$$\left(\int_{-A}^A \left|\frac{d}{d\xi}[\tilde{\varphi}_n(\xi) - e^{-\xi^2/2}]\right|^2 d\xi\right)^{1/2} \leq 0.76\Gamma_n^4\left(\int_{-A}^A \xi^6(1+2\xi^2+\xi^4)e^{-\xi^2}\, d\xi\right)^{1/2}$$

$$\leq 3.8\sqrt{\Gamma(7/2)}\,\Gamma_n^4. \qquad \square$$

**A.3.** We gathered here certain facts which are needed in order to complete the proof of Theorem 2.1.

REMARK A.1. The restriction to $U^c$ of the conditional characteristic function $\varphi^*$ belongs to $L^1(\mathbb{R})$.

Recalling the definition of the conditional distribution of (13), given $\beta$, one can write

$$\varphi^*(\xi) = \prod_{j=1}^{\nu} \varphi_0(\pi_j(\beta) \cdot \xi) = o\left(\frac{1}{\xi^{18}}\right),$$

where the latter equality, in view of the tail assumption (3), holds as $|\xi| \to +\infty$.

REMARK A.2. Assume that the real-valued random variable $Z$ satisfies $\mathsf{E}(Z^2) = 1$ and its characteristic function $\varphi$ belongs to $L^1(\mathbb{R})$; then $\varphi$ is an element of $H^1(\mathbb{R})$ and inequality

$$\int_a^{+\infty} (\varphi'(\xi))^2 \, d\xi \leq |\varphi(a)| + \int_a^{+\infty} |\varphi(\xi)| \, d\xi$$

holds for any $a$.

To prove this fact, take $b > a$ and integrate by parts to obtain

$$\int_a^b (\varphi'(\xi))^2 \, d\xi \leq |\varphi(a)\varphi'(a)| + |\varphi(b)\varphi'(b)| + \int_a^b |\varphi(\xi)| \cdot |\varphi''(\xi)| \, d\xi$$

$$\leq |\varphi(a)| + |\varphi(b)| + \int_a^b |\varphi(\xi)| \, d\xi.$$

To complete the argument, observe that the law of $Z$ is absolutely continuous [since $\varphi$ is in $L^1(\mathbb{R})$] which, in turn, implies that $|\varphi(b)| \to 0$ as $b \to +\infty$ in view of the Riemann–Lebesgue lemma.

PROOF OF (34). Notice that

$$\prod_{j=1}^{\nu}(1 + \pi_j^2 \eta^2) = \sum_{k=0}^{\nu} E_k(\pi_1^2, \ldots, \pi_\nu^2) \eta^{2k}$$

holds true in view of the definition of elementary symmetric function $E_k$. See, for example, Section 1.9 of Merris (2003). From the definition of $U$, it follows that $\nu$ is strictly greater than $\overline{n}$ on the complement of $U$. Hence,

$$\prod_{j=1}^{\nu}(1 + \pi_j^2 \eta^2) \geq E_{\overline{n}}(\pi_1^2, \ldots, \pi_\nu^2) \eta^{2\overline{n}},$$

and to complete the proof, it remains to show that $E_{\overline{n}} \geq \overline{\varepsilon}(\overline{n})$ is valid for any $\overline{n}$. With the aim of verifying such an inequality, recall (11) to write

$$
\begin{aligned}
1 = \left(\sum_{j=1}^{\nu} \pi_j^2\right)^m &= \sum_{i_1+\cdots+i_\nu=m} \frac{m!}{i_1!\cdots i_\nu!}(\pi_1^{2i_1}\cdots\pi_\nu^{2i_\nu}) \\
&\geq m! E_m(\pi_1^2,\ldots,\pi_\nu^2) \qquad (m \in \mathbb{N})
\end{aligned}
\tag{45}
$$

and use this fact to prove the desired claim by induction. It is trivially true for $\overline{n} = 1$. For $\overline{n} \geq 2$, the *Newton identities* with $M_h := \sum_{j=1}^{\nu} \pi_j^{2h}, h = 1, 2, \ldots,$ give

$$
\begin{aligned}
E_h &= \frac{1}{h}\left\{E_{h-1} + (-1)^{h-1}\left[M_h + \sum_{j=1}^{h-2}(-1)^j E_j M_{h-j}\right]\right\} \\
&\geq \frac{1}{h}\left\{\left(\frac{1}{(h-1)!} - 2^{h-2}\overline{\delta}\right) + (-1)^{h-1}\left[M_h + \sum_{j=1}^{h-2}(-1)^j E_j M_{h-j}\right]\right\}
\end{aligned}
$$

[by the inductive hypothesis]

$$
\begin{aligned}
&\geq \frac{1}{h}\left(\frac{1}{(h-1)!} - 2^{h-2}\overline{\delta}\right) - \overline{\delta} - \sum_{j=1}^{h-2}\overline{\delta} \qquad \text{[from (45)]} \\
&\geq \overline{\varepsilon}.
\end{aligned}
$$

To verify the last inequality for $h = 3, 4, \ldots,$ one can recall that $(2^{h-2}+h-1)$ is not greater than $2^{h-1}$. $\square$

**Acknowledgments.** The authors are grateful to Federico Bassetti for a few significant bibliographical suggestions. They wish to thank the reviewers for their time, comments and suggestions. They would also like to thank Eric Carlen and Maria Carvalho for their invaluable discussions.

E. Dolera
E. Gabetta
E. Regazzini
Dipartimento di Matematica "Felice Casorati"
Università degli Studi di Pavia
via Ferrata 1
27100 Pavia
Italy
E-mail: dolera.emanuele@libero.it
ester.gabetta@unipv.it
eugenio.regazzini@unipv.it